\documentclass[a4paper,11pt]{article}
\usepackage{amsmath,amsfonts,amsthm,amssymb}
\usepackage{graphicx,url,subfig}
\graphicspath{ {./Images/} }
\usepackage{color,tikz,mathtools}
\usepackage{float,wrapfig}
\usepackage[margin = 1.0in]{geometry}

\usepackage{thm-restate}

\newcommand{\R}{\mathbb{R}}

\newcommand{\N}{\mathbb{N}}
\newcommand{\F}{\mathcal{F}}

\newcommand{\ol}{\overline}
\newtheorem{theorem}{Theorem}

\newcounter{cases}
\newcounter{subcases}[cases]

\usepackage{color}

\usepackage[colorlinks=true,citecolor=blue]{hyperref}

\newcommand{\diam}{\mathrm{diam}}

\newcommand{\NN}{\mathbb{N}}

\newcommand{\cC}{\mathcal{C}}

\newcommand{\cF}{\mathcal{F}}

\newcommand{\cL}{\mathcal{L}}

\newcommand{\cO}{\mathcal{O}}

\newcommand{\cS}{\mathcal{S}}

\newcommand{\remove}[1]{}

\newtheorem{remark}[theorem]{Remark}

\title{
Helly theorem for affine spaces without dimensions
}

\author{
Sutanoya Chakraborty\footnote{Indian Statistical Institute, Kolkata, India}
\and
Arijit Ghosh\footnotemark[1]
\and
Soumi Nandi\footnote{
The Institute of Mathematical Sciences, Chennai, India
}
}

\date{}

\begin{document}
\maketitle

\begin{abstract}
    We prove a no-dimensional Helly theorem for affine spaces and convex sets using the unboundedness framework of Aronov, Goodman, and Pollack ({\em Computational Geometry}, 2002). This generalizes the fundamental result of Adiprasito, Bárány, Mustafa, and Terpai on the no-dimensional Helly theorem for points and convex sets ({\em Discrete \& Computational Geometry}, 2020). Additionally, we establish the optimality of our result. 
\end{abstract}
\section{Introduction}

Helly's theorem states that if every subset of size at most $d+1$ of a finite collection \( \mathcal{F} \) of convex sets in \( \mathbb{R}^d \) has a common point, then the entire collection $\mathcal{F}$ has a common point~\cite{Helly23}.
The finiteness of $\F$ can be removed if we assume that the convex sets in $\F$ are compact, therefore unless otherwise stated explicitly, we will not assume anything about the cardinality of a collection of convex sets.
Since its discovery, Helly's theorem has led to numerous applications and generalizations~\cite{Eckhoff1993helly,AmentaLP2017,DeLoeraGMM2019discrete,BaranyK2022helly}.

A set \( T \subseteq \mathbb{R}^d \) is said to \emph{pierce} a family \( \mathcal{F} \) of subsets of \( \mathbb{R}^d \) if every set in \( \mathcal{F} \) has a nonempty intersection with \( T \). 
A typical goal in a {\em Helly-type theorem} is to establish the existence of an integer $h(k, d)$ such that if any $h(k, d)$ sets in a collection $\mathcal{F}$ can be pierced by a $k$-flat, then the entire family also can be pierced by a \( k \)-flat. However, Santal\'{o}~\cite{Santal2009} proved that such a result is impossible even in the case of piercing convex sets with a line.

Hadwiger~\cite{Hadwiger56} showed that if $\mathcal{F}$ is a countable collection of disjoint convex sets in $\mathbb{R}^{d}$, each congruent to a fixed compact convex set with nonempty interior, and if every $d+1$ sets from $\mathcal{F}$ can be pierced by a line, then the entire family $\F$ can be pierced by a line. Danzer, Grünbaum, and Klee~\cite{DGK} later relaxed the congruence assumption, requiring the convex sets to be compact and have bounded diameters.

Further progress was made by Aronov, Goodman, Pollack, and Wenger~\cite{ABJP}, who established the first Helly-type result for hyperplanes, focusing on families of \emph{well-separated} compact convex sets in higher dimensions. Aronov, Goodman, and Pollack~\cite{AronovGP2002} extended these results to general $k$-flats for all $0 \leq k \leq d-1$, for collections of compact convex sets that are \emph{unbounded} in $k$-independent directions. Our approach builds on the framework and proof techniques developed by Aronov, Goodman, and Pollack.

Let $B(p, r)$ denote the \emph{closed ball} of radius $r \geq 0$ centered at $p \in \R^d$, and for \( S \subseteq \R^d \) and \( q \in \R^d \), the \emph{distance} from \( q \) to \( S \) is given by $d(q, S) := \inf_{x \in S} \|q - x\|$.
Adiprasito, B\'{a}r\'{a}ny, Mustafa, and Terpai~\cite{AdiprasitoBMT2020} introduced the first no-dimensional variant of the classical Helly’s theorem.

\begin{theorem}[Adiprasito, B\'{a}r\'{a}ny, Mustafa, and Terpai~\cite{AdiprasitoBMT2020}]\label{Point}
Let \(\cF\) be a collection of compact convex sets in \(\R^d\), $b \in \R^{d}$, and \(1 \leq r \leq d\). Suppose that for every choice of \(r\) sets \(C_1, C_2, \dots, C_r\) in \(\cF\), their intersection contains a point from $B(b,1)$. 
There exists a point \( q \in \mathbb{R}^d \) such that for all $C \in \mathcal{F}$, we have $d(q, C) \leq \dfrac{1}{\sqrt{r}}$.
\end{theorem}

They also proved the following colorful variant of the above theorem. 
\begin{theorem}[Adiprasito, B\'{a}r\'{a}ny, Mustafa, and Terpai~\cite{AdiprasitoBMT2020}]\label{PointColor}
	Let $\F_{1}, \F_{2}, \dots, \F_{r}$ be families of compact convex sets in $\R^d$ with $r\leq d$, and 
    $b \in \R^{d}$. 
	Suppose that for every $r$-tuple $\left(C_{1}, C_{2}, \dots,C_{r}\right) \in \F_{1} \times \F_{2} \times \dots \times \F_{r}$, we have 
    $\left(\bigcap\limits_{1 \leq i \leq r }C_{i}\right) \bigcap B(b,1) \neq \emptyset$.
    Then there exist a point $q\in\R^d$ and an index $i \in \left\{ 1, \dots, r\right\}$ such that, for all $C\in\cF_i$, we have $d(q,C)\leq \dfrac{1}{\sqrt{r}}$.
\end{theorem}

\begin{remark}\label{remark_no-dimensional_Helly_theorem}
    The original proofs of Theorem~\ref{Point} and Theorem~\ref{PointColor} apply to finite families of convex sets. However, standard point-set topology arguments extend them to infinite families of compact convex sets. 
\end{remark}

We generalize the above results for points to $k$-flats. We also prove some impossibility results which establish the optimality of our generalization.

\paragraph{Notations.}  
We use the following notations throughout the paper.
\begin{itemize}  
    \item 
        We denote the set of all real numbers and positive integers by $\R$ and $\N$ respectively. 

    \item The origin of $\mathbb{R}^{d}$ is denoted by $\cO$.  

    \item For \( n \in \mathbb{N} \), we define \( [n] := \{1, \dots, n\} \).  

    \item The {\em Euclidean distance} between \( p, q \in \mathbb{R}^{d} \) is denoted by \( \|p - q\| \).  

    \item 
        The {\em linear span} of a set $S \subseteq \R^{d}$ is defined as  
        $$
            \mathrm{span}(S) : = \left\{ \sum\limits_{i =1}^{n}\lambda_{i}x_{i} \mid n \in \N, x_{1}, x_{2}, \dots, x_{n} \in S, \lambda_{1}, \lambda_{2}, \dots, \lambda_{n} \in \R \right\}.
        $$

    \item For a set \( X \), \( |X| \) denotes its cardinality.  

    \item The {\em diameter} of a set \( C \subseteq \mathbb{R}^{d} \) is  
    $\mathrm{diam}(C) := \sup \left\{ \|c_1 - c_2\| \mid c_1, c_2 \in C \right\}$.

    \item 
        For \( b \in \mathbb{R}^{d} \) and \( r > 0 \), the {\em closed ball} and {\em open ball} of radius \( r \) centered at \( b \) are defined as $B(b, r) := \left\{ p \in \mathbb{R}^{d} \mid \|p - b\| \leq r \right\}$ and $B^{o}(b, r) := \left\{ p \in \mathbb{R}^{d} \mid \|p - b\| < r \right\}$, respectively.

    \item The {\em unit sphere} in \( \mathbb{R}^{d} \), centered at the origin, is denoted by \( \mathbb{S}^{d-1} \).  

    \item The {\em distance} between two subsets \( S_1, S_2 \subseteq \mathbb{R}^{d} \) is defined as 
    \[
        d(S_1, S_2) := \inf\left\{ \|p_1 - p_2\| \mid p_1 \in S_1, p_2 \in S_2 \right\}.
    \]  

    \item For \( p, q \in \mathbb{R}^{d} \), the closed line segment connecting \( p \) and \( q \) is denoted by \( \overline{pq} \).  
\end{itemize}

\section{Our results}\label{defth}

Before we give the statements of our results, we need to give some required definitions. 
We define the {\em central projection map} $f:\R^{d}\setminus\{\cO\} \to \mathbb{S}^{d-1}$,\remove{ where $\mathbb{S}^{d-1}$ is the unit sphere centered at $O$,} in the following way: for all $x \in \R^{d}\setminus\{\cO\}$, $f: x \mapsto {x}/{\|x\|}$.

We say $y \in \mathbb{S}^{d-1}$ is a {\em limiting direction} of a collection $\F$ of subsets of $\R^{d}$
if there exist two infinite sequences $\left\{ s_{n}\right\}_{n \in \mathbb{N}}$ and $\left\{ S_{n}\right\}_{n \in \N}$ satisfying the following properties:
\begin{itemize}
    \item 
        for all $n \in \N$, $S_{n} \in \F$

    \item 
%       \color{red} 
       for all $n \neq m$ in $\N$, $S_{n} \neq S_{m}$ 
       
%       \color{blue} Is this condition necessary? I think we should remove this \color{black}
    
    \item 
        for all $n \in \N$, $s_{n} \in S_{n}$
        
    \item
        $\lim\limits_{n \to \infty} \|s_{n}\| = \infty$ and 
        $\lim\limits_{n \to \infty} f(s_{n}) = y$
\end{itemize}
\noindent
For a collection $\F$ of subsets of $\R^{d}$, {\em limiting direction set} $LDS(\F)$ of $\F$ is defined as 
$$
    LDS(\F) : = \left\{ y \in \mathbb{S}^{d-1} \, : \, y\; \mbox{is a limiting direction of} \; \F \right\}.
$$
A collection $\F$ is {\em $k$-unbounded} if the vector space spanned by $LDS(\F)$ has dimension at least $k$. Aronov, Goodman, and Pollack~\cite{AronovGP2002} introduced this notion to prove a Helly theorem for $k$-flats, where $0 \leq k < d-1$. 
We establish the following generalization of their result.

%Aronov, Goodman and Pollack~\cite{AronovGP2002} proved the following Helly-type theorem for $k$-flats:
%\begin{theorem}[Aronov, Goodman and Pollack~\cite{AronovGP2002}]\label{kflat}
%	Let $\F$ be a family of compact convex sets in $\R^d$ such that there exists an $R\in\R^{+}$ for which $diam(C)<R$ $\forall C\in\F$ and $\cup_{C\in\F}C$ is unbounded. Now, 
%$\F$ is $k$-unbounded if the limiting direction set of $\F$ has dimension greater than or equal to $k$. In addition, 
%if $\F$ is $k$-unbounded and any $d+1$ sets in $\F$ can be intersected by a $k$-flat, then there exists a $k$-flat that intersects all sets in $\F$.
%\end{theorem}

%One can prove a {\em colorful version} of the above theorem.

\begin{restatable}[Colorful Helly's theorem for $k$-flats]{theorem}{colHellykflats}
    Let $\F_{1}, \F_{2}, \dots, \F_{d+1}$ be $k$-unbounded families of compact convex sets in $\R^d$ with $0 \leq k \leq d-1$, and there exists $R > 0$ such that for all $C \in \cF_{1} \cup \F_{2} \cup \dots \cup \cF_{d+1}$, $\diam(C) < R$. Also, assume that for every $(d+1)$-tuple $(C_{1}, C_{2}, \dots, C_{d+1}) \in \cF_{1}\times \F_{2} \times \dots \times \cF_{d+1}$ there exists a $k$-flat $K$ such that $K\cap C_{i} \neq \emptyset$, for all $i \in [d+1]$. Then, there exist an index $j\in[d+1]$ and a $k$-flat $\widetilde{K}$ such that for all $C \in \cS_{j}$, we have $\widetilde{K} \cap C \neq \emptyset$ .
 \label{col_aronov1}
 \end{restatable}
The following colorful theorem is the main technical contribution of this paper. The no-dimensional Helly theorem for $k$-flats follows as a direct consequence.

\begin{restatable}{theorem}{colnodimHellykflats}
    Let $\F_{1}, \F_{2}, \dots, \F_{r}$ be families of convex sets in $\R^{d}$, $k < r \leq d$, and the family $\F:= \F_{1} \cup \F_{2} \cup \dots \cup \F_{r}$ of compact convex sets satisfy the following properties:
    \begin{itemize}
        \item[(i)] 
            There exists $R > 0$ such that for all $C \in \F$ we have $\diam(C) < R$.

        \item[(ii)]
            There exist a set $J_{k} = \{ j_{1}, \dots, j_{k}\} \subseteq [r]$ and a collection of $k$ linearly independent vectors $\left\{ y_{j_{1}}, \dots, y_{j_{k}}\right\} \subseteq \mathbb{S}^{d-1}$ such that $y_{j_{i}} \in LDS\left( \F_{i} \right)$ for all $i \in [k]$.
    \end{itemize}
    Additionally, for any $r$-tuple $\left(C_{1}, C_{2}, \dots, C_{r}\right) \in \F_{1}\times \F_{2} \times \dots \times \F_{r}$ there exists a $k$-flat that intersects the closed unit ball $B\left(\cO,1\right)$ and every convex set $C_i$ for all $i \in [r]$. Then, there exist a $k$-flat $K$ and an index $j \in [r]$ such that for every $C \in \mathcal{F}_j$, we have 
    $d(C, K) \leq \sqrt{\dfrac{1}{r - k}}$.
    \label{main-result}
\end{restatable}
The following two results are immediate consequences of the above theorem.

\begin{theorem}[No-dimensional Helly theorem for $k$-flats]
    \label{noDimKflat}
    Let $\F$ be a $k$-unbounded family of compact convex sets in $\R^d$, and there exists an $R > 0$ such that for all $C \in \F$, we have $\diam(C) < R$.
    For $r \in \N$ with $k < r \leq d$ and $b \in \R^{d}$, if for every $C_{1}, C_{2}, \dots,C_r$ in $\F$ there exists a $k$-flat that intersects $B(\cO,1)$ and every $C_i$ for all $i\in\{1,2,\dots,r\}$, then there exists a $k$-flat $K$ such that, for all $C \in \F$, we have 
    $d(C,K) \leq \sqrt{\dfrac{1}{r-k}}$.
\end{theorem}

\begin{theorem}[Colorful no-dimensional Helly theorem for $k$-flats]
    \label{noDimKflat_colorful}
    Let $\F_{1}, \F_{2}, \dots, \F_{r}$ be $k$-unbounded families of compact convex sets in $\R^{d}$, $k < r \leq d$, and there exists $R >0$ such that for all $C \in \F_{1} \cup \F_{2} \cup \dots \F_{r}$, we have $\diam(C) < R$. Additionally, if for every $r$-tuple $(C_{1}, C_{2}, \dots, C_{r}) \in \F_{1} \times \F_{2} \times \dots \F_{r}$ there exists a $k$-flat that intersects $B(\cO,1)$ and every $C_i$ for all $i\in\{1,2,\dots,r\}$. Then, there exist a $k$-flat $K$ and an index $i \in [r]$ such that, for all $C \in \F$, we have 
    $d(C,K) \leq \sqrt{\dfrac{1}{r-k}}$.
\end{theorem}

\noindent 
% If each $\F_{i}$ in Theorem~\ref{main-result} is $k$-unbounded then we get the following colorful generalization of the above Theorem~\ref{noDimKflat}.

% \begin{theorem}[Colorful generalization of Theorem~\ref{noDimKflat}]\label{colorNoDimKflat}
%     Let $\F_{1}, \dots, \F_{r}$ be $k$-unbounded families of convex sets in $\R^d$ where $k<r\leq d$, and there exists $R > 0$ such that $\forall C \in \F_{1}\cup \dots \cup \F_{r}$ we have $\diam(C) < R$. 
%     If for any $r$-tuple $\left(C_{1}, \dots, C_{r}\right) \in \F_{1}\times \dots \times \F_{r}$ there exists a $k$-flat that intersects the closed unit ball $B\left(0,1\right)$ and every convex set $C_i$ for all $i \in \{1, \dots, r \}$, then there exists a $k$-flat $K$ and $j \in \{1, \dots, r\}$ such that, for all $C \in \F_{j}$, we have
%     \begin{equation*}
%         d(C,K) \leq \sqrt{\frac{1}{r-k}} \; .
%     \end{equation*}
% \end{theorem}

In the theorems above, we assume that the convex sets have bounded diameter. This condition is essential and cannot be relaxed. Consider hyperplanes in $\mathbb{R}^d$. Any finite collection of hyperplanes can be pierced by a line passing through the origin. However, for any $k$-flat $K$ with $k \leq d-1$ and any $\Delta > 0$, there exists a hyperplane $H$ such that $d(K, H) > \Delta$. The following two theorems complement the above no-dimensional Helly theorems by demonstrating the tightness of our bounds and the necessity of the $k$-unboundedness condition.

\begin{restatable}[On families not being $k$-unbounded]{theorem}{famnotkunbdd}
There exists a family $\cF$ of convex sets in $\R^{3}$ such that 
    \begin{itemize}
        \item 
            there exists $R > 0$ such that $\diam(C) < R$ for all $C \in \cF$,

        \item
            $\cF$ is $1$-unbounded,
            
        \item
            any three convex sets in $\cF$ can be pierced by a plane ($2$-dimensional affine space) passing through the origin $\cO$, and
            
        \item
            for any plane $K$ in $\R^{3}$ there exists a $C_{K} \in \cF$ such that $d\left( K, C_{K}\right) > 1$.
    \end{itemize}
\label{impossibility1}
\end{restatable}

\begin{restatable}[Tightness of the bound in Theorem~\ref{main-result}]{theorem}{tightness}
There exist families $\cF_{1}, \cF_{2}, \cF_{3}$ of compact convex sets in $\R^{3}$ such that 
    \begin{itemize}
        \item 
            for all $C \in \cF_{1} \cup \cF_{2} \cup \cF_{3}$, $\diam(C) = \sqrt{2}$,
                
        \item 
            both $\cF_{1}$ and $\cF_{2}$ are $1$-unbounded,
                
        \item
            for all $\left(C_{1}, C_{2}, C_{3} \right) \in \cF_{1} \times \cF_{2} \times \cF_{3}$ there exists a line $L$ that pierces $C_{1}, C_{2}, C_{3}$, and $d(L, \cO) \leq 1$, and
        
        \item
            for every line $K$ in $\R^{3}$ and for each $j \in [3]$, we have 
            $\max\limits_{C \in \cF_{j}} d(C, K) \geq \dfrac{1}{\sqrt{2}}$.
            
    \end{itemize}
\label{impossibility2}\end{restatable}

\section{Proofs of the claimed results}
\label{sec-proof}

\subsection{Helly theorems for $k$-flats: colorful and no-dimensional}
\label{ssec:helly_theorems_colorful_no-dimensional}

Lov\'{a}sz, and later B\'{a}r\'{a}ny~\cite{Barany82}, introduced the {\em colorful Helly theorem}, a fundamental generalization of Helly's theorem. This result is essential for completing the proof of Theorem~\ref{col_aronov1}.

\begin{theorem}[Lov\'{a}sz and B\'{a}r\'{a}ny~\cite{Barany82}: Colorful Helly Theorem]\label{th:barany}
    Let \(\mathcal{F}_1, \mathcal{F}_2, \dots, \mathcal{F}_{d+1}\) be families of compact convex sets in \(\mathbb{R}^d\). Suppose that for every \((d+1)\)-tuple \((C_1, \dots, C_{d+1}) \in \mathcal{F}_1 \times \dots \times \mathcal{F}_{d+1}\), the intersection \(\bigcap_{i=1}^{d+1} C_i\) is nonempty. Then, there exists an index \(j \in [d+1]\) such that \(\mathcal{F}_j\) can be pierced by a single point.
\end{theorem}

\begin{remark}\label{remark_colorful_Helly_theorem}
    Note that B\'{a}r\'{a}ny's proof of Theorem~\ref{th:barany}~\cite{Barany82} applies to finite families of convex sets in $\mathbb{R}^d$. As with Theorem~\ref{Point} and Theorem~\ref{PointColor}, standard arguments from point-set topology allow us to extend it to infinite families of compact convex sets.

\end{remark}

We first prove the following colorful generalization of the Helly theorem for $k$-flats proved by Aronov, Goodman and Pollack~\cite{AronovGP2002}.

\begin{proof}[Proof of Theorem~\ref{col_aronov1}]
Suppose for each $i\in[d+1],\; \cL_i$ is a set of $k$ linearly independent unit vectors in the limiting directions set $LDS(\cF_i)$ of $\cF_i$. Observe that there exists a linearly independent set of $k$ unit vectors, denoted by $\cL=\{z_1,\dots, z_k\}$, such that for each $i\in [k],\;z_i\in\cL_i$. Let $K$ denote the $k$-dimensional vector subspace $\mathrm{span}\left( \cL\right)$ of $\R^{d}$. Now for each $i\in [k]$, since $z_i\in\cL_i$, therefore there exist two sequences 
$\{C_{i,n}\}_{n\in\NN}$ and $\{ x_{i,n} \}_{n \in \NN}$ such that 
\begin{itemize}
    \item 
        for all $n \in \NN$, $C_{i,n} \in \cF_{i}$ and $x_{i,n} \in C_{i,n}$,
        
    \item 
        $\lim_{n\to \infty} \| x_{i,n} \| = \infty$ and $\lim_{n \to \infty} f(x_{i,n}) = z_{i}$.
\end{itemize}
Recall the definition of projection map $f: \R^{d}\setminus \{O\} \to \mathbb{S}^{d-1}$ with $f: x \longmapsto {x}/{\|x\|}$.

For each $i\in [d+1], i>k$, take a $C_i\in\cF_i$. Then $C_{k+1},\dots, C_{d+1}$ together with $C_{1,n},\dots, C_{k,n}$ is a colorful $(d+1)$-tuple and therefore pierceable by a $k$-flat, say $K_n$. Let $a_{n} \in K_{n}$ be a point in $K_{n}$ that is closest to $\cO$. Since the diameter of the sets $C_{k+1},\dots, C_{d+1}$ is bounded by $R$, there exists a $R' > 0$ such that 
$$
    \bigcup\limits_{k < i \leq d+1} C_{i} \subseteq B \left( \cO, R' \right),
$$ 
and therefore $\|a_{n}\| \leq R'$, for all $n \in \mathbb{N}$. As $\{ a_{n} \}_{n \in \mathbb{N}} \subseteq B(\cO, R')$, the sequence $\{ a_{n} \}_{n \in \mathbb{N}}$ has a limit point, say, $ a \in \R^{d}$. Without loss of generality, we can assume that $\lim_{n \to \infty} a_{n} = a$.

Observe that, for all $n \in \NN$, $K_{n} = a_{n} + K'_{n}$ where $K'_{n}$ is a $k$-dimensional subspace of $\R^{d}$ parallel to $K_{n}$. Now for each $i \in [k]$, let $y_{i,n} \in C_{i,n} \cap K_{n}$. For $i\in [k]$, the following infinite sequences 
$$
    \{ y_{i,n} - a_{n}\}_{n \in \N} \mbox{ and } \left\{ f(y_{i,n} - a_{n}) :=  \frac{y_{i,n} - a_{n}}{\| y_{i,n} - a_{n} \|} \right\}_{n \in \N}
$$
satisfy
\begin{itemize}
    \item {\bf Property-1:}
        for all $i \in [k]$, $\lim\limits_{n \to \infty} f(y_{i,n} - a_{n}) = z_{i}$, and 

    \item {\bf Property-2:}
        there exists $n_{0} \in \N$ such that for all $n \geq n_{0}$, we have
        $K_{n} = a_{n} + \mathrm{span} \left\{ y_{i, n} - a_{n}\, \mid \, i \in [k] \right\}$.

%    \item 
%        for all $i \in [k]$, $\lim\limits_{n \to \infty} f(y_{i,n} - a_{n}) = z_{i}$.
        
\end{itemize}
Note that {\bf Property-2} follows from {\bf Property-1}. We now establish {\bf Property-1}. Using the facts that $\|x_{i,n} - y_{i,n}\| \leq R$ (as the diameter of $C_{i,n}$ is at most $R$), $\|a_{n}\| \leq R'$, and $\lim_{n\to \infty} \|x_{i,n}\| = \lim_{n\to \infty} \|y_{i,n} - a_{n}\| = \infty$, observe that 
\begin{align*}
&\lim_{n \to \infty} \left\| \frac{y_{i,n} - a_{n}}{\|y_{i,n} - a_{n}\|} - \frac{x_{i,n}}{\|x_{i,n}\|}\right\| & \\
&= \lim_{n \to \infty} \left\| \frac{\|x_{i,n}\| y_{i,n} \ - \ \|x_{i,n}\|a_{n} \ - \ \|y_{i,n} - a_{n}\| x_{i,n}}{\|y_{i,n} -a_{n}\| \ \|x_{i,n}\|} \right\| & \\
&\leq \lim_{n\to \infty} \left\| \frac{\|x_{i,n}\| y_{i,n} \ - \ \|y_{i,n} - a_{n}\| x_{i,n}}{\|y_{i,n} -a_{n}\| \ \|x_{i,n}\|} \right\|
+ \lim_{n \to \infty} \frac{\|a_{n}\|}{\|y_{i,n} -  a_{n}\|} &\mbox{using triangle inequality}\\
&\leq \lim_{n\to \infty} \left\| \frac{\|x_{i,n}\| y_{i,n} \ - \ \|y_{i,n} - a_{n}\| y_{i,n} \ - \ \|y_{i,n} - a_{n}\| (x_{i,n} - y_{i,n})}{\|y_{i,n} -a_{n}\| \ \|x_{i,n}\|} \right\| &\mbox{as $\lim_{n\to \infty} \frac{\|a_{n}\|}{\|y_{i,n} - a_{n}\|} = 0$}\\
&\leq \lim_{n\to \infty} \left\| \frac{\|x_{i,n}\| y_{i,n} \ - \ \|y_{i,n} - a_{n}\| y_{i,n} \ - \ \|y_{i,n} - a_{n}\| (x_{i,n} - y_{i,n})}{\|y_{i,n} -a_{n}\| \ \|x_{i,n}\|} \right\| &\\
&\leq
\lim_{n\to \infty} \left\| \frac{\|x_{i,n}\| y_{i,n} \ - \ \|y_{i,n} - a_{n}\| y_{i,n} }{\|y_{i,n} -a_{n}\| \ \|x_{i,n}\|} \right\|
+
\lim_{n\to \infty} 
\frac{\|x_{i,n} - y_{i,n}\|}{\|x_{i,n}\|} & \mbox{using triangle inequality}\\
&\leq
\lim_{n\to \infty} \left\| \frac{\|x_{i,n}\| y_{i,n} \ - \ \|y_{i,n} - a_{n}\| y_{i,n} }{\|y_{i,n} -a_{n}\| \ \|x_{i,n}\|} \right\| &\mbox{as $\lim_{n\to \infty} \frac{\|x_{i,n} - y_{i,n}\|}{\|x_{i,n}\|} = 0$}\\
&\leq
\lim_{n\to \infty} \frac{ \|y_{i,n}\| \  \left|\Big(\|x_{i,n}\| - \|y_{i,n} - a_{n}\| \Big)\right| }{\|y_{i,n} -a_{n}\| \ \|x_{i,n}\|} &\\
&\leq
\lim_{n\to \infty} \frac{ \Big(\|y_{i,n} - a_{n}\| + \|a_{n}\| \Big) \  \|x_{i,n} - y_{i,n} + a_{n}\|}{\|y_{i,n} -a_{n}\| \ \|x_{i,n}\|} &\mbox{using triangle inequality}\\
&\leq 
\lim_{n\to \infty} \frac{ \Big(\|y_{i,n} - a_{n}\| + \|a_{n}\| \Big) \  \Big( \|x_{i,n} - y_{i,n}\| + \|a_{n}\| \Big)}{\|y_{i,n} -a_{n}\| \ \|x_{i,n}\|} &\mbox{using triangle inequality}\\
&\leq
\lim_{n \to \infty} \left( 1 + \frac{\|a_{n}\|}{\|y_{i,n} - a_{n}\|} \right) \times \lim_{n\to \infty}\left( \frac{\|x_{i,n} - y_{i,n} \|}{\|x_{i,n}\|} + \frac{\|a_{n}\|}{\|x_{i,n}\|} \right) &\\
&= 0,
\end{align*}
the last equality follows from the fact that 
$$
\lim_{n \to \infty} \frac{\|a_{n}\|}{\|y_{i,n} - a_{n}\|} = 
\lim_{n\to\infty}\frac{\|a_{n}\|}{\|x_{i,n}\|} = \lim_{n\to \infty} \frac{\|x_{i,n} - y_{i,n}\|}{\|x_{i,n}\|} = 0.
$$
Therefore, as $n \to \infty$, $f(y_{i,n} - a_{n}) = \frac{y_{i,n} - a_{n}}{\|y_{i,n} - a_{n}\|}$ converge to $z_{i} = \lim_{n\to \infty}f(x_{i,n}) = \lim_{n\to \infty}\frac{x_{i,n}}{\|x_{i,n}\|}$. This completes the proof of {\bf Property-1}.

Using the above facts, we get that the following sequence $\{ K_{n} \}_{n \in \N}$ of $k$-flats convergences to the $k$-flat $a + K$. Consider the following sequences $\left\{ C_{i} \cap K_{n} \right\}_{n \in \N}$, for all $ k < i \leq d+1$, of subsets $C_{i}$. Since $C_{i}$ is a compact set and $C_{i} \cap K_{n}$ are non-empty closed subsets of $C_{i}$ therefore $C_{i} \cap \left( a + K\right) \neq \emptyset$. Therefore, $a+K$ pierces $C_{k+1}, \dots, C_{d+1}$.

%\color{blue}
%We assume that $K_n=a_n+K'_n$, where $K'_n$ is a $k$-flat passing through the origin and $|a_n|=d(K_n,O)$. Since each $S_{i,n}$ is of bounded diameter containing $x_{i,n}$ such that $f(x_{i,n})\rightarrow z_i$, as $n\rightarrow \infty$, where $\|z_i\|=1$, $a_n$ must be bounded, for all $n\in\NN$. Then the sequence $\{a_n\}_{n\in\NN}$ must have a limit point, say, $a\in\RR^d$.
%Now for each $i\in [k]$, suppose $y_{i,n}\in S_{i,n}\cap K_n$. Then $y_{i,n}-a_n\in K'_n$ and also $f(y_{i,n}-a_n)=\frac{y_{i,n}-a_n}{\|y_{i,n}-a_n\|}\in K'_n$. Again since each $S_{i,n}$ is of bounded diameter containing $x_{i,n}$ such that $\|x_{i,n}\|\rightarrow\infty$, as $n\rightarrow\infty$, we have $\|y_{i,n}\|\rightarrow\infty$, as $n\rightarrow\infty$. Then $f(y_{i,n}-a_n)=\frac{y_{i,n}-a_n}{\|y_{i,n}-a_n\|}\rightarrow z_i-a$, as $n\rightarrow\infty$. Then by compactness of the sets $C_{k+1},\dots,C_{d+1}$, we get that $C_{k+1},\dots,C_{d+1}$ is pierceable by the $k$-flat $a+\cK$.

% This follows that $B_{k+1},\dots, B_{d+1}$ can be pierced by a $k$-flat arbitrarily close to the direction of $\cK$. So by compactness of $B_i$'s we can say that, $B_{k+1},\dots, B_{d+1}$ can be pierced by a $k$-flat in the direction of $\cK$, i.e, parallel to $\cK$.

%\color{black}
Suppose that, for each $i \in [d+1]$, the family $\mathcal{F}'_{i}$ is obtained by orthogonally projecting the sets of $\mathcal{F}_i$ onto the $(d-k)$-dimensional subspace $K^{\perp}$, which is the orthogonal complement of $K$.
Every colorful $(d-k+1)$-tuple from $\cF'_{k+1}, \dots, \cF'_{d+1}$ is pierceable by a point in the space $K^{\perp}$. 
Using colorful Helly theorem (Theorem~\ref{th:barany}), there exists a $i\in \{ k+1, \dots, d+1\}$ such that $\cF'_i$ is pierceable by a point in $K^{\perp}$. Hence there exists a $k$-flat parallel to $K$ that pierces all the convex sets in the family $\cF_i$.
\end{proof}

We now give the proof of the generalization of the no-dimensional Helly theorem of Adiprasito, B\'{a}r\'{a}ny, Mustafa, and Terpai~\cite{AdiprasitoBMT2020} for $k$-flats. 

\begin{proof}[Proof of Theorem~\ref{main-result}]
%Given a set $F\subseteq \RR^d$, $\overline{F}$ denotes the closure of $F$. Since $d(p, F) = d(p, \overline{F})$ for any point $p \in \mathbb{R}^d$, it suffices to show that there exist an index $j \in [r]$ and a $k$-flat $K$ such that
%$$
%    d(\overline{C}, K) \leq \sqrt{\frac{1}{r-k}}
%$$
%for every \( C \in \mathcal{F}_j \).
%
%        
Without loss of generality assume that $J_k=[k]$. For each $i \in [k]$, there exist two infinite sequences $\left\{ F_{i, n} \right\}_{n \in \N}$ and $\left\{ x_{i,n} \right\}_{n \in \N}$ such that
\begin{itemize}
    \item 
        for all $n \in \NN$, $F_{i,n} \in \cF_{i}$ and $x_{i,n} \in F_{i,n}$,
        
    \item 
        $\lim_{n\to \infty} \| x_{i,n} \| = \infty$ and $\lim_{n \to \infty} f(x_{i,n}) = y_{i}$.  
\end{itemize}
Let $K$ denote the subspace $\mathrm{span}\left(\{y_1,\dots,y_k\}\right)$. For each $i\in [r],\; i>k,$ take any $F_i\in\cF_i$. Observe that for all $n\in\NN$, there exists a $k$-flat $K_n$ intersecting $B(\cO,1)$ and piercing $F_{1,n},F_{2,n},\dots, F_{k,n}, F_{k+1},\dots, F_r$. 
Using the same arguments as in the proof of Theorem~\ref{col_aronov1}, we can show that there exists a point $a$ in $B(\cO,1)$ such that the convex sets $F_{k+1}, F_{k+2}, \dots, F_{r}$ is pierced by the $k$-flat $a + K$.

%This implies that $F_{k+1},\dots, F_r$ can be pierced by a $k$-flat arbitrarily close to the direction of $K$. Now by compactness of $\Bar{F_i}$'s, we conclude that $\exists a\in K^{\perp}$ such that the $k$-flat $a+K$ pierces $\Bar{F}_{k+1},\dots,\Bar{F}_r$.Since $a+K$ intersects $B(O,1)$, we must have $\|a\|\leq 1$.

Let $\cF'_{i}$ denote the family obtained by orthogonally projecting the sets in $\cF_{i}$ onto the $(d-k)$-dimensional subspace $K^{\perp}$. Observe that for all $\left( C'_{k+1}, C'_{k+2}, \dots, C'_{r}\right) \in \cF'_{k+1} \times \cF'_{k+2} \times \dots \times \cF'_{r}$, we have 
$$
    \left(\bigcap_{i=k+1}^{r} {C'}_{i} \right)\bigcap B( \cO,1 )\neq\emptyset.
$$
By Theorem~\ref{PointColor}, there exist a point $q\in K^{\perp}$ and an index $i\in\{k+1,\dots, r\}$ such that for all $C' \in\cF_i$, $d \left( q, \pi \left( C' \right) \right) \leq \sqrt{\dfrac{1}{r-k}}$. Therefore, for all $C \in \cF_{i}$, we have $d (C, K+q) \leq \sqrt{\dfrac{1}{r-k}}$.
\end{proof}

\subsection{Necessity of $k$-unboundedness in Theorem~\ref{impossibility1}}

%\famnotkunbdd*

%\begin{figure}
%    \centering
%    \includegraphics[scale = 0.80]{old_version/newfigure1.pdf}
%    \caption{}
%    \label{fig-convex-set-being-bounded}
%\end{figure}

\begin{figure}
    \centering
    \includegraphics[scale = 1.250]{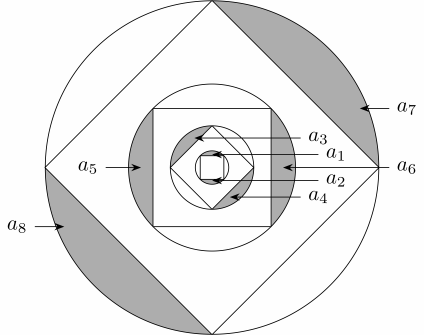}
    \quad
    \includegraphics[width = 6.0cm, height= 6.0cm]
    {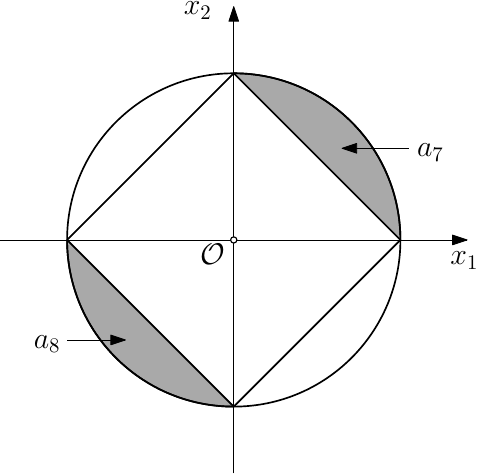}
    \caption{An example demonstrating the necessity of $k$-unboundedness condition. This figure has been taken from~\cite{AronovGP2002}.}
    \label{fig:k1}
\end{figure}

To establish the necessity of $k$-unboundedness in Theorem~\ref{impossibility1}, we build on a construction given by Aronov, Goodman, and Pollack~\cite{AronovGP2002}.

\begin{proof}[Proof of Theorem~\ref{impossibility1}]
Consider the eight shaded convex regions in Figure~\ref{fig:k1}, formed by four circles and four squares centered at a point $\cO$. Let $x_3 = 0$ be the plane containing Figure~\ref{fig:k1}, and assume without loss of generality that $\cO$ is the origin in $\mathbb{R}^{3}$. We label these eight shaded convex regions as $a_1, a_2, \dots, a_8$. Observe that any three of these convex sets can be intersected by a straight line passing through $\cO$.

We construct additional convex sets by arranging the eight given sets in a fixed order and sequentially elevating them along the $x_{3}$-axis in integer steps, that is, the additional convex sets will be of the form $a_{i} + (0,0,8n + i -1)$, with $i \in [8]$ and $n \in \N$. 
This process ensures that for any $n \in \mathbb{N}$, infinitely many sets lie outside $B(\cO, n)$. As a result, we obtain a countably infinite sequence $\mathcal{F}$ of sets, where $\mathcal{F}$ is $1$-unbounded but not $2$-unbounded.
To see why this is true, consider any infinite sequence of points \( \{p_n\}_{n \in \mathbb{N}} \) chosen from the sets in \( \mathcal{F} \). Since we are elevating eight bounded sets in the \( x_1x_2 \)-plane to the planes \( x_3 = n \) for all \( n \in \mathbb{N} \), the points \( p_n \) take the form \( p_n = (x_1, x_2, n) \). 
Since \( |x_1| \) and \( |x_2| \) are bounded, both  
$\frac{x_1}{\sqrt{x_1^2 + x_2^2 + n^2}}$ and $\frac{x_2}{\sqrt{x_1^2 + x_2^2 + n^2}}$ tend to zero as \( n \to \infty \). Thus,  
$f(p_n) := \frac{p_n}{\| p_n \|} \to (0, 0, 1)$ as $n \to \infty$.

    Let $\ell$ denote the length of the side of the smallest square in Figure \ref{fig:k1}.
    We show that the distance of any plane is greater than $1$ from at least one set in $\F$ when $\ell$ is large enough.
    Let $K$ be a plane for which the maximum distance from the sets in $\F$ is minimized, and $l_{K}$ be the intersection of $K$ with the plane $x_{3} = 0$. Clearly, $K$ must be perpendicular to the plane that contains the first $8$ sets, because otherwise for any $R>0$ we would find a set $C$ in $\F$ for which $d(K,C)>R$.
    Consider the straight line $l_{K}$ that is the intersection of $K$ and the plane given by the equation $x_{3} = 0$.
    If $l_{K}$ is moved on the $x_{3} = 0$ plane closer to $\cO$ along the line perpendicular to $l_{K}$ from $\cO$, the quantity $\max\left\{d\left(a_{2n},l_{K}\right),d\left(a_{2n-1},l_{K}\right)\right\}$ does not increase for $n \in \left\{ 1,2,3,4 \right\}$.
    Since $d(K,a_i)=d(l_{K},a_i)$ for all $i\in[8]$, we can take $K$ to be passing through $\cO$. Let the side-lengths of the $4$ squares in Figure \ref{fig:k1} be $\ell_{1}\, (=\ell), \ell_{3}, \ell_{5}, \ell_{7}$, where the side of side-length $\ell_{i}$ is shared by the set $a_{i}$.
    Let the diagonals of the largest square in Figure \ref{fig:k1} lie on the $x_{1}$-axis and the $x_{2}$-axis respectively.
    Then, if $l_{K}$ makes an angle $\theta\in[0,\pi)$ with the $x_{1}$-axis, then we have the following: for $i \in \left\{ 1, 3, 5, 7 \right\}$ we have 
    \begin{align}
        d\left(a_{i}, l_{K} \right) = d \left( a_{i+1}, l_{K} \right),
    \end{align}
    and
    \begin{align}
        d(a_{1},l_{K}) &=
                \begin{cases}
                    \frac{\ell_1}{\sqrt{2}}\sin\left({\pi}/{4}-\theta \right) & \text{if } 0\leq\theta\leq {\pi}/{4}\\
                    \frac{\ell_1}{\sqrt{2}} \sin\left(\theta-3\pi/4\right) & \text{if } {3\pi}/{4} \leq \theta \leq \pi \\
                    0 & \text{otherwise}
                \end{cases}
    \end{align}
    \begin{align}
        d(a_{3},l_{K}) &=
                \begin{cases}
                    \frac{\ell_3}{\sqrt{2}}\sin\theta & \text{if } 0 \leq \theta \leq {\pi}/{4} \\
                    \frac{\ell_3}{\sqrt{2}}\sin\left({\pi}/{2}-\theta\right) & \text{if } {\pi}/{4} \leq \theta \leq {\pi}/{2} \\
                    0 & \text{otherwise}
                \end{cases}
    \end{align}
    \begin{align}
        d(a_{5},l_{K}) &=
                \begin{cases}
                    \frac{\ell_5}{\sqrt{2}}\sin\left(\theta - {\pi}/{4}\right) & \text{if } {\pi}/{4} \leq \theta \leq {\pi}/{2}\\
                    \frac{\ell_5}{\sqrt{2}}\sin \left( {3\pi}/{4}-\theta \right) & \text{if }{\pi}/{2} \leq \theta \leq {3\pi}/{4}\\
                    0 & \text{otherwise}
                \end{cases}
    \end{align}
    \begin{align}
        d(a_7,l_{K}) &=
                \begin{cases}
                    \frac{\ell_7}{\sqrt{2}}\sin\left(\theta - {\pi}/{2}\right) & \text{if }{\pi}/{2} \leq \theta \leq {3\pi}/{4} \\
                    \frac{\ell_7}{\sqrt{2}} \sin \left(\pi-\theta\right) & \text{if } {3\pi}/{4} \leq \theta \leq \pi \\
                    0 & \text{otherwise}
                \end{cases}
    \end{align}
    Now, observe that 
    \begin{align*}
        \sum_{i=1}^{8} d\left(a_{i},l_{K}\right) \geq \frac{\ell}{\sqrt{2}} \times 
        \min_{\theta\in \left[0, {\pi}/{4}\right]} \left(\sin\theta + \sin\left({\pi}/{4}-\theta \right)\right)
        \geq \ell\sqrt{2}\sin\left({\pi}/{8}\right)
    \end{align*}
    Thus, for $\ell>\sqrt{2}/\sin\left({\pi}/{8}\right)$, there are no planes that are at most $1$ distance away from each set in $\F$.
\end{proof}

% \begin{theorem}[Tightness of the bound in Theorem~\ref{main-result}: Restatement of Theorem~\ref{impossibility2}]\label{restatement-impossibility2}
% There exist families $\cF_{1}, \cF_{2}, \cF_{3}$ of convex sets in $\R^{3}$ satisfying the following properties:
%     \begin{itemize}
%         \item 
%             $\forall C \in \cF_{1} \cup \cF_{2} \cup \cF_{3}$, $\diam(C) = \sqrt{2}$,
                
%         \item 
%             both $\cF_{1}$ and $\cF_{2}$ are $1$-unbounded,
                
%         \item
%             $\forall (C_{1}, C_{2}, C_{3}) \in \cF_{1} \times \cF_{2} \times \cF_{3}$ there exists a line $L$ that pierces $C_{1}, C_{2}, C_{3}$, and $d(L, \mathcal{O}) \leq 1$ and
        
%         \item
%             for every line $K$ in $\R^{3}$ and for every $j \in [3]$, 
%             $\max\limits_{C \in \cF_{j}} d(C, K) \geq \frac{1}{\sqrt{2}}$.
%     \end{itemize}
% \end{theorem}

\subsection{Tightness of Theorem~\ref{main-result}}

%\tightness*

\begin{figure}
\centering
\includegraphics[width = 5.0cm, height = 9.50cm]{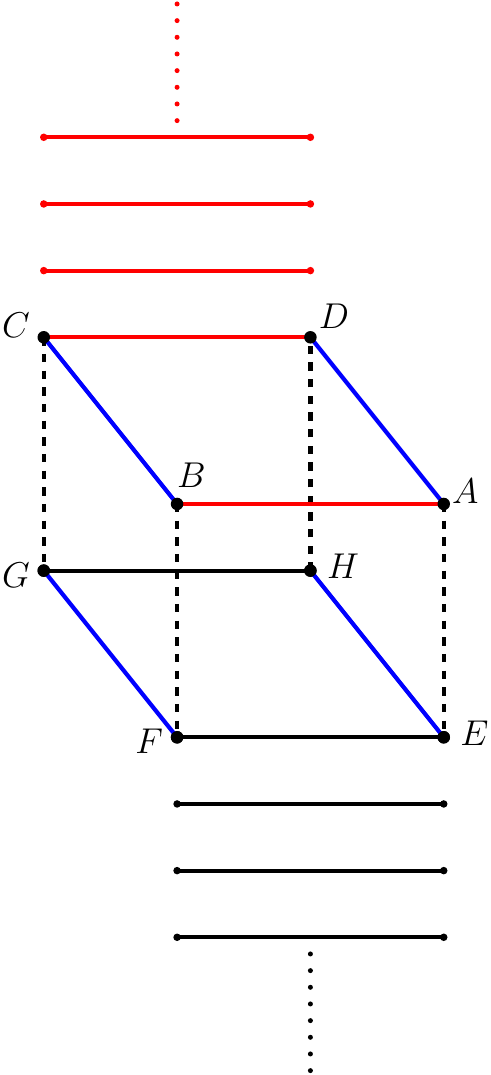}
\caption{An example demonstrating the tightness of the bound given in Theorem~\ref{impossibility2}.}
\label{fig:bound2}
\end{figure}

We provide an explicit extremal example establishing the tightness of the bound 
given in Theorem~\ref{main-result}.

\begin{proof}[Proof of Theorem~\ref{impossibility2}]    
Let $\{A,B,C,D,E,F,G,H\}$ be the $8$ vertices cube in $\R^3$ whose centroid is the origin $\cO$ and side length is $\sqrt{2}$ and $EFGH$ is parallel to the plane $x_3=0$ (see Figure~\ref{fig:bound2}).
Define $\F_{i} := \{S_{i,1},S_{i,2},S_{i,3},\dots\}$ for $i \in \{1,2,3$\} in the following way: Set $S_{1,1}=\ol{AB}$, $S_{1,2}=\ol{CD}$ and for $n>2$, $S_{1,n}$ is the line segment $\ol{CD}$ raised to the height $x_3=n$. Similarly, set $S_{2,1}=\ol{GH}$, $S_{1,2}=\ol{EF}$ and for $n>2$, $S_{2,n}$ is the line segment $\ol{EF}$ lowered to the height $x_3=-n$. Now set $S_{3,1}=\ol{BC},S_{3,2}=\ol{DA},S_{3,3}=\ol{FG},S_{3,4}=\ol{HE}$, and let $S_{3,4n+j}$ be the set $S_{3,j}$, for all $j \in \{0,1,2,3\}$ (see Figure~\ref{fig:bound2}). Clearly, any colorful $3$-tuple $(C_1, C_2, C_3)$, $C_i\in\F_i$, $i\in[3]$, can be hit by a straight line that is at most at a distance $1$ away from the centroid $\cO$.
To see this, consider a colorful $3$-tuple $\cC=(C_1, C_2, C_3)$, $C_i\in\F_i$, $i\in[3]$, and a line $\ell$ piercing the $3$ sets in $\cC$. Then $\ell$ must intersect two opposite sides of a facet of the cube $ABCDEFGH$ with centroid $\cO$ and side length $\sqrt{2}$. Therefore, the distance of the origin $\cO$ from that line cannot be greater than the distance of the origin from one of the edges of the facet, which is $1$.

Now for each colorful $3$-tuple $\cC=(C_1,C_2,C_3)$, $C_i\in\F_i$, $i\in[3]$, suppose $\ell_{\cC}$ is the line transversal of $\cC$ closest to the centroid $\cO$ and $\cL=\left\{\ell_{\cC}\;|\;\cC\text{ is a colorful $3$-tuple }\right\}$.\remove{Let $\mathcal{L}$ denote the set of all straight line transversals $l$ of colorful $3$-tuples such that $l$ is as close to $\cO$ as possible.} Let $l_n$ denote the straight line transversal in $\mathcal{L}$ that passes through $S_{1,n},S_{2,1},S_{3,1}$. Clearly, $d(l_n,\cO)\to 1$ as $n\to\infty$.
Then 
$$
    \sup_{l\in\mathcal{L}}d(\cO,l) = 1.
$$
Now note that for the straight line $L$ that is perpendicular to the plane on which $ABCD$ lies and passes through $\cO$, we have $$\inf_{i\in[3]}\sup_{C\in\F_i}d(l,C)\ge\inf_{i\in[3]}\sup_{C\in\F_i}d(L,C),$$ for any straight line $l$ in $\R^3$. We can show this in the following way: let $l_1$ be a straight line such that $$\sup_{C\in\F_1}d(l_1,C)=\inf_{l\in\mathcal{L}}\sup_{C\in\F_1}d(C,l).$$
Then $l_1$ must be perpendicular to the plane on which $ABCD$ lies, otherwise, the supremum of its distances from sets in $\F_1$ would be infinity. Then, $l_1$ must be equidistant from both $S_{1,1}$ and $S_{1,2}$, and therefore, we can take $l_1$ to be $L$. Similar arguments show that 
$$
    \sup_{C\in\F_2}d(L,C)=\inf_{l\in\mathcal{L}}\sup_{C\in\F_1}d(C,l).
$$ 
We have $$\inf_{i\in[3]}\sup_{C\in\F_i}d(L,C)=\frac{1}{\sqrt{2}}.$$
To see that $$\inf_{l\in\mathcal{L}}\sup_{C\in\F_3}d(C,l)=\frac{1}{\sqrt{2}},$$ project $\ol{BC}$, $\ol{AD}$, $\ol{FG}$, $\ol{EH}$ onto the plane $P$ that contains $ABEF$. If there is a straight line $l_3$ such that 
$$
    \sup_{C\in\F_3}d(l_3,C)=\inf_{l\in\mathcal{L}}
    \sup_{C\in\F_3}d(C,l),
$$ 
then let the projection of $l_3$ onto $P$ be $l'_3$. If $L$ is not the straight line that minimizes $$\inf_{l\in\mathcal{L}}\sup_{C\in\F_3}d(C,l),$$ then the perpendicular distance from $l'_3$ to $A,B,E$ and $F$ must be smaller than $\frac{1}{\sqrt{2}}$. Let, without loss of generality, $A$ be the point from which $l'_3$ is the farthest. Then, we must have another point among $B,E,$ and $F$ from which $l'_3$ has the same distance as $A$. This point then must be $E$, because otherwise, we could have taken $L$ to be $l_3$. This means that $l'_3$ passes through the centroid of the square $ABEF$ and two points from $A,B,E,$ and $F$ lie on each side of $l'_3$. But this implies that $l'_3$ must be parallel to $AB$ since $l'_3$ has the minimum distance from both $A$ and $B$ and is at least as close to $E$ and $F$, which is a contradiction. Therefore 
$$
    \inf_{i\in[3]}\sup_{C\in\F_i}d(L,C)=\frac{1}{\sqrt{2}}.
$$ 
\end{proof}

\bibliographystyle{alpha}
\bibliography{reference}

\begin{thebibliography}{DLGMM19}

\bibitem[ABMT20]{AdiprasitoBMT2020}
K.~A. Adiprasito, I.~B{\'{a}}r{\'{a}}ny, N.~H. Mustafa, and T.~Terpai.
\newblock {Theorems of Carath{\'{e}}odory, Helly, and Tverberg Without
  Dimension}.
\newblock {\em Discrete \& Computational Geometry}, 64(2):233 -- 258, 2020.

\bibitem[ADLS17]{AmentaLP2017}
N.~Amenta, J.~A. De~Loera, and P.~Sober{\'o}n.
\newblock {Helly’s Theorem: New Variations and Applications}.
\newblock {\em {Algebraic and Geometric Methods in Discrete Mathematics}},
  685:55--95, 2017.

\bibitem[AGP02]{AronovGP2002}
B.~Aronov, J.~E. Goodman, and R.~Pollack.
\newblock {A Helly-type Theorem for Higher-dimensional Transversals}.
\newblock {\em Computational Geometry}, 21(3):177--183, 2002.

\bibitem[AGPW00]{ABJP}
B.~Aronov, J.~E. Goodman, R.~Pollack, and R.~Wenger.
\newblock {A Helly-Type Theorem for Hyperplane Transversals to Well-Separated
  Convex Sets}.
\newblock In {\em Proceedings of the 16th Annual Symposium on Computational
  Geometry, SoCG}, page 57–63, 2000.

\bibitem[B{\'{a}}r82]{Barany82}
I.~B{\'{a}}r{\'{a}}ny.
\newblock {A Generalization of Carath{\'{e}}odory's Theorem}.
\newblock {\em Discrete Mathematics}, 40(2-3):141--152, 1982.

\bibitem[BK22]{BaranyK2022helly}
I.~B{\'a}r{\'a}ny and G.~Kalai.
\newblock {Helly-type Problems}.
\newblock {\em Bulletin of the American Mathematical Society}, 59(4):471--502,
  2022.

\bibitem[DGK63]{DGK}
L.~Danzer, B.~Gr{\"u}nbaum, and V.~Klee.
\newblock {\em Helly's Theorem and Its Relatives}.
\newblock Proceedings of Symposia in Pure Mathematics: Convexity. American
  Mathematical Society, 1963.

\bibitem[DLGMM19]{DeLoeraGMM2019discrete}
J.~De~Loera, X.~Goaoc, F.~Meunier, and N.~Mustafa.
\newblock {The discrete yet ubiquitous theorems of Carath{\'e}odory, Helly,
  Sperner, Tucker, and Tverberg}.
\newblock {\em Bulletin of the American Mathematical Society}, 56(3):415--511,
  2019.

\bibitem[Eck93]{Eckhoff1993helly}
J.~Eckhoff.
\newblock {Helly, Radon, and Carath{\'e}odory Type Theorems}.
\newblock In {\em Handbook of Convex Geometry}, pages 389 -- 448. Elsevier,
  1993.

\bibitem[Had56]{Hadwiger56}
H.~Hadwiger.
\newblock {\"{U}ber einen Satz Hellyscher Art}.
\newblock {\em Arch. Math.}, 7:377--379, 1956.

\bibitem[Hel23]{Helly23}
E.~Helly.
\newblock {\"{U}ber Mengen konvexer K\"{o}rper mit gemeinschaftlichen Punkten}.
\newblock {\em Jahresbericht der Deutschen Mathematiker-Vereinigung}, 32:175 --
  176, 1923.

\bibitem[San40]{Santal2009}
L~Santal{\'o}.
\newblock Un teorema sobre conjuntos de paralelepipedos de aristas paralelas.
\newblock {\em Publ. Inst. Mat. Univ. Nac. Litoral}, 2:49--60, 1940.

\end{thebibliography}

\appendix

\end{document}